\documentclass[12pt]{amsart}
\usepackage{amsfonts}
\usepackage{amscd}
\usepackage{amsmath, amsfonts,amssymb,amsthm,mathrsfs,xypic}
\usepackage{graphicx}
\xyoption{curve}
\usepackage{fancybox}
\newtheorem{thm}{Theorem}[section]

\newtheorem{lem}[thm]{Lemma}

\newtheorem{prop}[thm]{Proposition}

\numberwithin{equation}{section}

\newcommand{\s}{\hfill\blacksquare}
\newcommand{\op}{{\operatorname{op}}}%%%
\newcommand{\End}{\operatorname{End}}%%%
\newcommand{\add}{\operatorname{add}}%%%
\newcommand{\Ker}{\operatorname{Ker}}%%%
\newcommand{\Hom}{\operatorname{Hom}}%%%
\newcommand{\Ext}{\operatorname{Ext}}%%%
\newcommand{\mo}{\operatorname{mod}}%%%
\newcommand{\proj}{\operatorname{proj}}%%%

\begin{document}

\title[Objective functors]{Objective triangle functors}
\author [Claus Michael Ringel and Pu Zhang] {Claus Michael Ringel and Pu Zhang}
\thanks{{\it 2010 Mathematical Subject Classification. \ 16E30, 18A22, 16E35.}}
\thanks{Supported by the NSF of China (11271251) and Specialized Research Fund for the Doctoral Program of Higher
Education (20120073110058).}

%%%%%%%%%%%%%%%%%%%%%ABSTRACT%%%%%%%%%%%%%%%%%%%%%%%%%%%%%%%%
\begin{abstract}
An (additive) functor $F:  \mathcal A \longrightarrow \mathcal B$ between additive categories 
is said to be objective, provided any morphism $f$ in $\mathcal A$ with $F(f) = 0$ factors through
an object $K$ with $F(K) = 0$.
In this paper we concentrate on triangle functors between
triangulated categories.  
The first aim of this paper is to characterize objective triangle
functors $F$ in several ways. Second, we are interested in the corresponding
Verdier quotient functors $V_F$, 
in particular we want do know under what
conditions $V_F$ is full. The third question to be considered concerns 
the possibility to factorize a given triangle functor $F = F_2F_1$
with $F_1$ a full and dense triangle functor and $F_2$ a faithful triangle functor.
It turns our that the behaviour of splitting monomorphisms (and splitting epimorphisms)
plays a decisive role.
\vskip5pt

{\it Key words and phrases: \   triangulated category, triangle
functor, objective functor, Verdier functor.}
\end{abstract}

\maketitle

%%%%%%%%%%%%%%%%%%%%%%%%%%%INTRODUCTION%%%%%%%%%%%%%%%%%%%%%%%%%%%%%%
\section{\bf Introduction}

Let $F:  \mathcal A \longrightarrow \mathcal B$ be a
functor between additive categories (all functors considered in this paper are
supposed to be covariant and additive). 
 Following [RZ], we say that $F$ is
 {\it objective}, provided any morphism $f: \ X
 \longrightarrow Y$ in $\mathcal A$ with $F(f) = 0$ factors through
 an object $K$ with $F(K) = 0$. We say that $F$ is {\it sincere},
 provided that $F$ sends non-zero objects to non-zero objects.
 Clearly, a functor is faithful if and only if it is
 objective and sincere.
 \medskip
 
In this paper we concentrate on triangle functors between
triangulated categories.  We will see that triangle functors behave
quite different from general (additive) functors between additive
categories, and also from exact functors between abelian
categories. For examples, a full functor between additive
categories may be not objective, but a full triangle functor between triangulated
categories is objective (see \ref{full}); on the other hand,
an exact functor between abelian categories is clearly objective (see
\ref{exactobj}), whereas
there are sincere triangle functors
between triangulated categories which are not objective (see section
\ref{sincere-notobjective}).
	\bigskip 

If $\mathcal I$ is an ideal of an additive category $\mathcal A$, then we denote by
$\mathcal A/\mathcal I$ the corresponding factor category: it has the same
objects as $\mathcal A$ and $\Hom_{\mathcal A/\mathcal I}(X,Y) =
\Hom_{\mathcal A}(X,Y)/\mathcal I(X,Y)$ for any pair $X,Y$ of objects in $\mathcal A$.
We denote by $\pi_{\mathcal I}: \mathcal A \to \mathcal A/\mathcal I$
the canonical projection functor, it sends an object $X$ to itself, and a
morphism $f$ to its residue class modulo $\mathcal I.$ 
Given a full subcategory
$\mathcal U$ of $\mathcal A$, we denote by $\langle\mathcal U \rangle$
the ideal generated by $\mathcal U.$ 

Assume now that $\mathcal A$ is a triangulated category and that $\mathcal K$ is
a triangulated subcategory of $\mathcal A$ (triangulated subcategories are always
assumed to be full subcategories). Then there is a 
triangulated category $\mathcal A/\mathcal K$ and a dense triangle functor
$V_{\mathcal K}:\mathcal A \longrightarrow \mathcal A/\mathcal K$
with the following universal property: 
$V_\mathcal K(\mathcal K) = 0$, and if $G:
\mathcal A \longrightarrow \mathcal B$ is a triangle functor with
$G(\mathcal K) = 0$, then there is a unique triangle functor $G': 
\mathcal A/\mathcal K \longrightarrow \mathcal B$ such
that $G = G' V_\mathcal K$ (see Verdier [V], or also Neeman [N]).
We call $V_{\mathcal K}$ the {\it Verdier quotient functor} for $\mathcal K$.
(There is no need to worry about a possible confusion using the same notation
$\mathcal  A/\mathcal  I$ and $\mathcal  A/\mathcal  K$, for ideals $\mathcal  I$ and triangulated subcategories
$\mathcal  K$, since a subcategory $\mathcal  U$ of $\mathcal  A$ is an ideal only in case $\mathcal  U = \mathcal  A$, see \ref{ideal}). 

If $F:\mathcal A \to 
\mathcal B$ is a functor between additive categories, then we
denote by $\ker(F)$ the class of morphisms $f$ in $\mathcal A$ such that
$F(f) = 0,$ and we denote by $\Ker(F)$ the full subcategory of $\mathcal A$ given by all
objects $X$ in $\mathcal A$ such that $F(X) = 0.$ 
Note that $\ker(F)$ is an ideal of $\mathcal A$, whereas $\Ker(F)$ is
a subcategory. Thus, given a functor $F$, we have two ideals $\langle \Ker(F)\rangle \subseteq
\ker(F).$ It is easy to see (\ref{objective}) that a functor $F$ is objective 
if and only if the ideals $\ker(F)$ and $\langle \Ker(F)\rangle$ coincide.
We will consider the factor category $\mathcal A/\ker(F)$ and we write
$\pi_F$ instead of $\pi_{\ker(F)}.$ 
If $F: \mathcal A\longrightarrow \mathcal B$ is a triangle
functor between triangulated categories, then the subcategory 
${\Ker}F$ 
is a triangulated subcategory of $\mathcal A$, thus we may consider the
Verdier quotient functor $V_{\Ker F}$, we will denote it by $V_F$.
Since $F({\Ker}F) = 0$, the
universal property of the Verdier quotient functor asserts that
there exists a unique triangle functor $\widetilde{F}: \mathcal
A/{\Ker} F \longrightarrow \mathcal B$, such that $F = \widetilde FV_F$,
the functor $\widetilde F$ is always sincere.
    \bigskip 

    The first aim of this paper is to characterize objective triangle
functors $F$ in several ways. Second, we are interested in the corresponding
Verdier quotient functors $V_F$, in particular we want do know under what
conditions $V_F$ is full. The third question to be considered concerns 
the possibility to factorize a given triangle functor $F = F_2F_1$
with $F_1$ a full and dense triangle functor and $F_2$ a faithful triangle functor.

Two conditions for a triangle functor $F$ will play a decisive role, namely 
the weak splitting monomorphism condition  ${\rm (WSM)}$ and
the isomorphism condition ${\rm (I)}$. If $F$ 
is faithful or full, then both conditions are satisfied (see \ref{faithful-RSM},
\ref{RSMimpliesSM} and \ref{WSMISM}).

${\rm (WSM)}$ \ \ For each morphism $u: X \longrightarrow Y$ in
$\mathcal A$ such that $F(u)$ is a splitting monomorphism in
$\mathcal B$, there exists a morphism $u': Y \longrightarrow X'$
such that $F(u'u)$ is an isomorphism in $\mathcal B$.

${\rm (I)}$ \ \ For each morphism $u: X \longrightarrow Y$ in
$\mathcal A$ such that $F(u)$ is an isomorphism in $\mathcal B$,
there exists a morphism $u': Y \longrightarrow X$ such that
$F(u)^{-1} = F(u').$

It is easy to see (see \ref{WSMISM}) that a functor $F$ satisfies both conditions
(WSM) and (I) if and only if it satisfies the splitting monomorphism condition (SM):

${\rm (SM)}$ \ \ For each morphism $u: X \longrightarrow Y$ in
$\mathcal A$ such that $F(u)$ is a splitting monomorphism in
$\mathcal B$, there exists a morphism $u': Y \longrightarrow X$ such
that $F(u'u) = 1_{F(X)}$.

Here are the main results of the paper:

%%%%%%%%%%%%%%%%%%%%%%Theorem 1.1%%%%%%%%%%%%%%%%%%%%%%%%%%%%%
\begin{thm} \label{1.1} \ Let $F: \mathcal
A \longrightarrow \mathcal B$ be a triangle functor between
triangulated categories. Then the following are equivalent:

$\rm (i)$ \ \ $F$ satisfies the condition ${\rm (WSM)}$;

$\rm (ii)$ \ \ $F$ is objective;

$\rm (iii)$ \ \ the induced functor $\widetilde{F}: \mathcal
A/{\Ker}F \longrightarrow \mathcal B$ is faithful.
\end{thm}

We say that an additive category $\mathcal A$ is a {\it Fitting category}
provided for any endomorphism $a:X \to X$ in $\mathcal A$ there exists 
a direct decomposition $X = X'\oplus X''$ with $a(X') \subseteq X', 
a(X'') \subseteq X''$ such that the restriction of $a$ to $X'$ is an automorphism
and the restriction of $a$ to $X''$ is nilpotent. For example, if 
$\mathcal A$ is a $\Hom$-finite $k$-category, where $k$ is a field, and any object of $\mathcal A$
is a finite direct sum of objects with local endomorphism rings, then $\mathcal A$ is a Fitting category (see
\ref{Fitting}).

%%%%%%%%%%%%%%%%%%%%%%Theorem 1.2%%%%%%%%%%%%%%%%%%%%%%%%%%%%%
\begin{thm} \label{1.2} \ Let $F: \mathcal A \to \mathcal B$ be a
triangle functor between triangulated categories.

$(1)$ \ \ If $V_F$ is full, then $F$ satisfies the condition ${\rm
(I)}$.

$(2)$ \ \  Assume that $F$ is objective or that
$\mathcal A$ is a Fitting category.
Then $V_F$ is full if and only if $F$ satisfies ${\rm (I)}$.
\end{thm}

%%%%%%%%%%%%%%%%%%%%%%Theorem 1.3%%%%%%%%%%%%%%%%%%%%%%%%%%%%%
\begin{thm} \label{1.3} \  Let $F: \mathcal A \longrightarrow \mathcal B$ be a
triangle functor between triangulated categories. Then the following
are equivalent:

$\rm (i)$ \ \  $F$ satisfies the condition ${\rm (SM)}$;

$\rm (ii)$ \ \ $F$ is objective and $V_F$ is full.

$\rm (iii)$ \ \ There is an equivalence of additive categories $\Phi:\mathcal A/\ker(F)
\to \mathcal A/\Ker(F)$ such that $V_F =\Phi \pi_F.$ 

$\rm (iv)$ \ \ There is factorization $F = F_2F_1$ where
$F_1$ is a full and dense triangle functor and $F_2$ is a faithful triangle functor.

$\rm (v)$ \ \ There is factorization $F = F_2F_1$ where
$F_1$ is a full triangle functor and $F_2$ is a faithful triangle functor.

\end{thm}

%%%%%%%%%%%%%%%%%%%%%%PRELIMINARIES%%%%%%%%%%%%%%%%%%%%%%%%%%
\section{\bf Preliminaries}

\noindent
{\bf Ideals and subcategories of additive categories.}

\begin{lem} \label{ideal} 
A subcategory $\mathcal U$ of an additive category $\mathcal A$
is an ideal if and only if $\mathcal U = \mathcal A$.
\end{lem} 

\noindent 
{\bf Proof.} Let $\mathcal U$ be an ideal of $\mathcal A$. 
If $X$ is an object in $\mathcal A$, then the zero map $0_X$ belongs to any ideal,
thus to $\mathcal U$. But since $\mathcal U$ is a subcategory, with $0_X$ also $1_X$ 
belongs to $\mathcal U$. The ideal generated by all the identity maps $1_X$ is clearly
$\mathcal A.$ $\s$
	  \bigskip 

\noindent
{\bf Objective functors.}
Let $F:\mathcal A \to \mathcal B$ be a functor
between additive categories. 
Recall that we denote by $\ker(F)$ the class of morphisms $f$ in $\mathcal A$ such that
$F(f) = 0$ (this is an ideal of the category $\mathcal A$) and by $\Ker(F)$ 
the full subcategory of $\mathcal A$ given by all
objects $X$ in $\mathcal A$ such that $F(X) = 0.$ Clearly, 
$\langle \Ker(F)\rangle \subseteq \ker(F).$ 

\begin{lem}\label{objective}
A functor $F$ is objective if and only if the ideals $\ker(F)$
and $\langle \Ker(F)\rangle$ coincide.
\end{lem}

\noindent 
{\bf Proof.} First, assume that $F$ is objective. Let $f:X \to Y$ belong to $\ker(F)$, thus 
$F(f) = 0.$ Since $F$ is objective, $f = hg$, with $g:X \to K,$ $h:K \to Y$ and
$F(K) = 0.$ Thus $K$ belongs to $\Ker(F)$ and therefore $f = hg = h\cdot 1_K\cdot g$
belongs to $\langle \Ker(F)\rangle$. 

Conversely, assume that $\ker(F) = \langle \Ker(F)\rangle$. Let 
$f:X \to Y$ be a morphism with $F(f) = 0$, thus $f$ belongs to $\ker(F)$ and
therefore to $\langle \Ker(F)\rangle$. This means that $f$ is of the form 
$\sum_{i=1}^m h_if_ig_i$ with maps $g_i:X \to K_i,$
$f_i:K_i \to K'_i,$ $h_i:K_i' \to Y$, where $K_i,K'_i$ are objects in 
$\Ker(F)$. Let $K = \bigoplus_{i=1}^m K_i$ and define maps 
$g = [g_1,\dots,g_m]^t:X \to K$ and $h = [h_1f_1,\dots,h_mf_m]:K \to Y$.
Then $f = hg$ shows that $f$ factors though the object $K$.
Of course,  $F(K) = 0.$ $\s$
   \medskip

\noindent
{\bf Triangulated categories and triangle functors.}
A triangulated category is of the form $\mathcal T = (\mathcal T, [1], \mathcal E)$
where $\mathcal T$ is an additive category, $[1]$ an automorphism of $\mathcal T$
and $\mathcal E$ a class of sixtuples of the form 
$X\stackrel u
\longrightarrow Y \stackrel v \longrightarrow Z 
\stackrel w\longrightarrow X[1]$ with objects $X,Y,Z$ and morphisms
$u,v,w$ (usually, we will denote such a sixtuple
just by $(X,Y,Z,u,v,w)$) satisfying some well-known axioms. 
The sixtuples in $\mathcal E$ are said to be the {\it distinguished triangles}. 
If $\mathcal A$  and $\mathcal B$ are triangulated
categories, a triangle functor from $\mathcal A$  to $\mathcal B$ is
a pair $F = (F, \xi)$, where $F: \mathcal A\longrightarrow \mathcal B$
is an additive functor, and $\xi: F\circ [1] \longrightarrow
[1]\circ F$ is a natural isomorphism, such that if 
$(X,Y,Z,u,v,w)$ is a distinguished triangle in $\mathcal A$, then
$(F(X),F(Y),F(Z),F(u),F(v),\xi_XF(w))$ is a
distinguished triangle in $\mathcal B$.  We should stress that given a triangle functor
$(F,\xi)$, there may not exist a triangle functor $(F',\xi')$ with $\xi'$ the
identity transformation such that $F$ and $F'$ are equivalent
(as pointed
out by Keller, see the appendix of Bocklandt [B]).
Note that 
triangle functors are also called exact functors or triangulated functors
(see e.g. [GM], [H]),[KS], [N]), but we follow the terminology used for example by 
Keller [Ke].
       \medskip 

\noindent
{\bf The Verdier quotient functor.}
Let us recall some
property of the Verdier quotient functor 
$V_\mathcal K:\mathcal A \to \mathcal A/\mathcal K$, namely that any morphism $x$ in $\mathcal A/\mathcal K$ 
can be written in the form $x = a/s = V_{\mathcal K}(a)V_{\mathcal K}(s)^{-1}$ 
where $a:X' \to Y$ and $s:X' \to X$ 
are morphisms in $\mathcal A$ such that there exists a distinguished triangle 
$(X',X,K,s,v,w)$ in $\mathcal A$ such that $K$ belongs to $\mathcal K$ (and also in the
form $x = V_{\mathcal K}(s')^{-1}V_{\mathcal K}(a')$ for some morphisms
$a':X\to Y'$ and $s':Y \to Y'$ in $\mathcal A$ with a distinguished triangle
$(Y,Y',K',s',v',w')$ in $\mathcal A$ such that $K'$ belongs to $\mathcal K$).
This follows directly from the construction of $\mathcal A/\mathcal K$ using the
calculus of fractions. 
	\medskip 

\noindent
{\bf Isomorphisms and splitting monomorphisms in triangulated categories.}
Recall that in a distinguished triangle $(X,Y,Z,u,v,w)$ the morphism
$u$ is a splitting
monomorphism if and only if $v$ is a splitting epimorphism, and if
and only if $w = 0$. Also, $u$ is an isomorphism if and only if
$Z = 0$. See Happel [H], I.1.4 and I.1.7.

%%%%%%%%%%%%%%%%%%%%CONDITIONS FOR TRIANGLE FUNCTORS%%%%%%%%%%%%%%%%%%%%%%%%%%%%%
\section{\bf Some conditions for triangle functors}

\noindent
Let $F: \mathcal A\longrightarrow \mathcal
B$ be a functor between additive categories.

\begin{prop} \label{WSMISM} Let $F: \mathcal A\longrightarrow \mathcal
B$ be a functor between additive categories.
Then F satisfies the conditions {\rm(WSM)} and {\rm (I)} if and only if
it satisfies the condition {\rm(SM)}. 
\end{prop}

\noindent {\bf Proof.}  First, assume that $F$ satisfies the condition (SM).
Let $u:X \to Y$ be a morphism in $\mathcal A$. If $F(u)$ is a splitting monomorphism,
then (SM) asserts the existence of $u'$ in $\mathcal A$ such that
$F(u'u) = 1_{F(X)}$, thus $F(u'u)$ is an isomorphism. This shows (WSM). 
If $F(u)$ is an isomorphism, then $F(u)$ is a splitting monomorphism, thus
there is $u$ in $\mathcal A$ such
that $F(u'u) = 1_{F(X)}$. Thus $F(u')F(u) = 1_{F(X)}$ and therefore
$F(u)^{-1} = F(u').$ This shows (I). 

Conversely, assume that the conditions (WSM) and (I) are satisfied. 
Let $F(u)$ be a splitting monomorphism with $u:
X\longrightarrow Y$. By ${\rm (WSM)}$ there exists a morphism $a: Y
\longrightarrow X'$ such that $F(au) = F(a)F(u)$ is an isomorphism
in $\mathcal B$. By ${\rm (I)}$ there exists a morphism $b: X'
\longrightarrow X$ such that $F(b)F(au) = 1_{F(X)}$. Put $u':
= ba: Y\longrightarrow X.$ Then $F(u')F(u) = F(b)F(au) = 1_{F(X)}.$ 
This proves that $F$ satisfies ${\rm (SM)}$. $\s$
     \bigskip
     
We need two further conditions for a functor $F:\mathcal A \to \mathcal B$.
   \medskip
   
${\rm (RSM)}$ \ \ {\it $F$ reflects splitting monomorphisms} (this means: if $u$
is a morphism in $\mathcal A$ such that
$F(u)$ is a splitting monomorphism in $\mathcal B$, then $u$ is a
splitting monomorphism in $\mathcal A$).
	  \medskip 

${\rm (RI)}$ \ \ {\it $F$ reflects isomorphisms} (this means: if $u$
is a morphism in $\mathcal A$ such that
$F(u)$ is an isomorphism in $\mathcal B$, then $u$ is an isomorphism in $\mathcal A$).
       \medskip 

\begin{prop} \label{faithful-RSM} Let $F: \mathcal A\longrightarrow \mathcal
B$ be a triangle functor between triangulated categories.
Then $F$ is faithful if and only if $F$ satisfies the conditions {\rm(RSM)}.
\end{prop}

\noindent {\bf Proof.} 
First, assume that $F = (F,\xi)$ is faithful. 
Let $u:X\to Y$ be a morphism in $\mathcal A$ such that $F(u)$ is a splitting
monomorphism. Let $(X,Y,Z,u,v,w)$ be a distinguished triangle in $\mathcal A$.
Then $(F(X),F(Y),F(Z),F(u),F(v),\xi_XF(w))$ is a distinguished triangle
in $\mathcal B$. Since $F(u)$ is a splitting monomorphism, $\xi_XF(w) = 0$,
thus also $F(w) = 0.$ Since $F$ is faithful, $w = 0$, thus $u$ is a splitting
monomorphism. Thus, the condition (RSM) is satisfied.

Conversely, assume that (RSM) holds. Let $w:Z \to X[1]$ be a morphism such that
$F(w) = 0.$ Then there is a  distinguished triangle $(X,Y,Z,u,v,w)$ in $\mathcal A$
and $(F(X),F(Y),F(Z),F(u),F(v),\xi_XF(w))$ is a distinguished triangle
in $\mathcal B$. Since $F(w) = 0$, also $\xi_XF(w) = 0$, thus $F(u)$ is 
a splitting monomorphism. Since $F$ satisfies the condition (RSM), 
$u$ is a splitting monomorphism,
thus $w = 0.$ This shows that $F$ is faithful.
$\s$ 
     \bigskip

\begin{prop} \label{RSMimpliesSM} Let $F: \mathcal A\longrightarrow \mathcal
B$ be a functor between additive categories. If $F$ is full or faithful, 
then it satisfies the condition {\rm (SM)}.
\end{prop}

\noindent {\bf Proof.} Let $u:X\to Y$ be a morphism in $\mathcal A$ such that
$F(u)$ is a splitting monomorphism, thus there is $b:F(Y) \to F(X)$ 
such that $bF(u) = 1_X$.  First, assume that $F$ is full. Then $b = F(u')$ for some 
$u':Y \to X$, and $F(u')F(u) = 1_{F(X)}$ shows that (SM) is satisfied.

Second, assume that $F$ is faithful, thus according to \ref{faithful-RSM},
$F$ satisfies the condition (RSM).
Let $u:X \to Y$ be a morphism in $\mathcal A$ such that $F(u)$ is a splitting
monomorphism. Since $F$ satisfies (RSM), 
there is $u'$ in $\mathcal A$ such that
$u'u= 1_X.$ Thus $F(u'u) = 1_{F(X)}.$   $\s$

\begin{prop} \label{RIequivalentsincere} Let $F: \mathcal A\longrightarrow \mathcal
B$ be a triangle functor between triangulated categories. Then $F$ is sincere if
and only if $F$ satisfies the condition {\rm (RI)}.
\end{prop}

\noindent 
{\bf Proof.} First, assume that $F$ is sincere. Let $u:X\to Y$ be a morphism in $\mathcal A$
such that $F(u)$ is an isomorphism. Let $(X,Y,Z,u,v,w)$ be a distinguished
triangle in $\mathcal A$. Then $(F(X),F(Y),F(Z),F(u),F(v),\xi_XF(w)$ is a distinguished
triangle in $\mathcal B$. Since $F(u)$ is an isomorphism, we have $F(Z) = 0$.
Since $F$ is sincere, this implies that $Z = 0$, thus $u$ is an isomorphism.
This shows that $F$ reflects isomorphisms.

Conversely, let us assume that $F$ reflects isomorphisms. In order to show that $F$
is sincere, let $Z$ be an object in $\mathcal A$ such that $F(Z) = 0.$ Consider
the map $u:Z \to 0$ in $\mathcal A$. If we apply $F$, we obtain a map
$F(u):F(Z) \to 0.$ Since $F(Z) = 0$, the map $F(u)$ is an isomorphism. Since
$F$ reflects isomorphisms, we see that $u$ itself is an isomorphism, but this
means that $Z = 0.$ $\s$

%%%%%%%%%%%%%%%%%%%%%Proof of Theorem 1.1%%%%%%%%%%%%%%%%%%%%%%%%%%%%%%%%
\section{\bf Objectivity of triangle functors.}

\noindent
The aim of this section is to prove Theorem 1.1 and to draw the attention to some
consequences. 

\subsection{The proof of Theorem 1.1.}

(i) $\implies$ (ii). We assume now that $F$ satisfies the condition (WSM).
We want to show that $F$ is objective, thus let $w$ be a morphism in $\mathcal A$
with $F(w) = 0,$ say $w:Z \to X[1]$. We take a distinguished triangle
$(X,Y,Z,u,v,w)$ in $\mathcal A$. Under $F$ we obtain the
distinguished triangle $(F(X),F(Y),F(Z),F(u),F(v),\xi_XF(w))$. Since $F(w) = 0$,
also $\xi_XF(w) = 0$, thus $F(u)$ is a split monomorphism. The condition
(WSM) provides a morphism $u':Y \to X'$ such that $F(u'u)$ is an isomorphism.
We need a distinguished triangle involving $u'u$, say $(X,X',K,u'u,f,h)$.
The given factorization of $u'u$ yields the following  commutative square
on the left:
$$ 
 \xymatrix{X \ar[r]^-u \ar @{=}[d] &
           Y \ar[d]^-{u'}\ar[r]^-{v} &
           Z \ar@{-->}[d]^-{g} \ar[r]^-w&
           X[1] \ar @{=}[d] 
           \\
           X \ar[r]^-{u'u} &
           X' \ar[r]^-{f}&
           K\ar[r]^h &
           X[1].}
$$ 
thus we obtain a morphism $g:Z \to K$ such that $w = hg.$ If we apply $F$
to the distinguished triangle $(X,X',K,u'u,f,h)$, we obtain the
distinguished triangle $(F(X),F(X'),F(K),F(u'u),F(f),\xi_XF(h))$. Since
$F(u'u)$ is an isomorphism, $F(K) = 0.$ Thus $w = hg$ is the required factorization.
	 \medskip

(ii) $\implies$ (iii). Assume that $F$ is objective. We want to show that $\widetilde F$
is faithful, thus consider a morphism $a/s$ in $\mathcal A/\Ker F$ with
$\widetilde F(a/s) = 0$. Here we deal with morphisms $s:X'\to X$ and $a:X' \to Y$ 
in $\mathcal A$ such that $V_F(s)$ is invertible. As a consequence,
also $F(s)$ is invertible.
We have 
$$
 F(a)F(s)^{-1} = \widetilde FV_F(a) \widetilde FV_F(s)^{-1} 
 = \widetilde F(V_F(a)V_F(s)^{-1}) 
 = \widetilde F(a/s) = 0,
$$ 
thus $F(a) = 0$. Since $F$ is objective, there is a factorization
$a = hg$, say with $g:X' \to K$ and $h:K \to Y$ such that $F(K) = 0.$ But
$F(K) = 0$ implies that $V_F(K) = 0$. 
Since $V_F(a) = V_F(h)V_F(g)$ factors through 
$V_F(K) = 0$, it follows that $V_F(a) = 0$, therefore also 
$a/s = V_F(a)V_F(s)^{-1} = 0.$ 
     \medskip

(iii) $\implies$ (i). 
We assume that $\widetilde F$ is faithful, thus objective.
Let $u:X\to Y$ 
be a morphism in $\mathcal A$ such that $F(u)$ is a splitting monomorphism.
Let $(X,Y,Z,u,v,w)$ be a distinguished triangle in $\mathcal A$, thus
$(F(X),F(Y),F(Z),F(u),F(v),\xi_XF(w))$ is a distinguished triangle in $\mathcal B$.
Since $F(u)$ is a splitting monomorphism, we know that $\xi_XF(w) = 0$, thus
$F(w) = 0$. Since $F = \widetilde FV_F$ and $\widetilde F$ is faithful, we see
that $V_F(w) = 0.$ It follows that $V_F(u)$ is a splitting monomorphism, thus there
is some $x$ in $\mathcal A/\Ker F$ such that $xV_F(u) = 1_{V_F(X)}.$ 
As we know, the morphism $x$
can be written in the form $x = V_F(s)^{-1}V_F(u')$ for some morphisms
$u':Y\to X'$ and $s:X \to X'$ in $\mathcal A$ with $V_F(s)$ invertible. 
This implies that $V_F(u'u) = V_F(u')V_F(u) = V_F(s)$ is an isomorphism.
If we now apply $\widetilde F$, we see that also $F(u'u) = \widetilde FV_F(u'u) =
\widetilde FV_F(s)$ is an isomorphism. $\s$
	   \medskip 

%%%%%%%%%%%%%%%%%%%%%%%%%%%%%%%%%%%%%%%%%%%%%%%%%%%%%%%%%%%%%
\subsection{The Verdier quotient functors are objective.}\label{Verdierfunctorisobjective}
As an immediate consequence of Theorem \ref{1.1}, we recover the following well-known result
(see, for example, Krause [Kr], Proposition 4.6.2):
      \medskip

\noindent
{\bf Corollary.}
{\it Let $\mathcal A$ be a triangulated category and $\mathcal K$
a triangulated subcategory of $\mathcal A$.  Then the  Verdier quotient
functor $V_\mathcal K: \mathcal A \longrightarrow \mathcal
A/\mathcal K$ is objective.}
	   \medskip 

\noindent 	   
{\bf Proof.} Clearly, $\widetilde V_{\mathcal K}$ is the identity functor on $\mathcal A/\mathcal K$, thus
faithful. It follows that $\widetilde V_{\mathcal K}$ is objective. $\s$

\subsection{Sincere triangle functors.} 
It is well-known that a sincere triangle functor $F$ which is full is also
faithful, see J. Rickard [Ric], p.446, l.1. The previous discussions provide
necessary and sufficient conditions for a sincere functor to be faithful:
Namely, for any triangle functor $F$, there are the following implications:
$$
 \text{faithful} \iff \text{(RSM)} \implies 
 \text{(SM)} \implies \text{(WSM)} \iff \text{objective} 
$$
(see \ref{faithful-RSM}, \ref{RSMimpliesSM}, \ref{WSMISM}, \ref{1.1}).
Since a sincere objective functor is of course faithful, all these conditions
are equivalent in case $F$ is sincere.

\subsection{Full triangle functors are objective.}
\label{full}\quad 
		  \medskip

{\bf Corollary.} {\it A full triangle functor is objective.}
     \medskip 

\noindent 
{\bf Proof.} According to Proposition \ref{RSMimpliesSM}, a full triangle functor 
satisfies the condition (SM), thus (WSM), and therefore $F$ is objective
by Theorem \ref{1.1}. $\s$

%%%%%%%%%%%%%%%%%%%%%Proof of Theorem 1.2, Part 1%%%%%%%%%%%%%%%%%%%%%%%%%%%%%%%%
\section{\bf Triangle functors $F$ with $V_F$ full}

The aim of this section is to present the proof of Theorem \ref{1.2}.

\subsection{\bf If $V_F$ is full, then (I) is satisfied}\label{}
We assume that $F:\mathcal A \to \mathcal B$
is a triangle functor such that $V_F$ is full. We want to show that $F$
satisfies the condition (I). 
Let  $u: X \to  Y$ be a morphism in $\mathcal A$
such that $F(u)$ is an isomorphism,  and $(X,Y,Z,u,v,w)$ a distinguished
triangle in $\mathcal A$. Thus \newline
$(F(X),F(Y),F(Z),F(u),F(V),\xi_XF(w))$
is a distinguished triangle in $\mathcal B$. Since $F(u)$ is an isomorphism,
$F(Z) = 0,$ thus $Z$ belongs to $\Ker(F)$ and therefore 
$V_F(u)$ is invertible in $\mathcal A/\Ker(F)$. Since $V_F$ is full, there is
a map $u':Y \to X$ such that $V_F(u') = (V_F(u))^{-1}$. If we apply the functor 
$\widetilde F$ (with $F = \widetilde FV_F$) we see that $F(u') = F(u)^{-1}$.

\subsection
{\bf The Verdier quotient functor $V_F$ for a functor $F$ satisfying (I)}\label{Verdier}
We assume that $F:\mathcal A \to \mathcal B$ is a triangle functor which satisfies the
condition (I). Note that 
the morphisms of $\mathcal A/\Ker(F)$ are of the form
$a/s =V_F(a)(V_F(s))^{-1}$,  where $a,s$ are morphisms in $\mathcal A$ 
with $V_F(s)$ being invertible. In order to show that $V_F$ is full, 
it is sufficient to show that the morphisms of the form $(V_F(s))^{-1}$ are
in the image of $V_F$. 

\subsection{\bf The case when $F$ is objective} 
Let $F:\mathcal A \to \mathcal B$ be an objective triangle functor which
satisfies the condition (I). Let $s$ be a morphism in $\mathcal A$ such that
$V_F(s)$ is invertible.  Apply the functor $\widetilde F$ (with
$\widetilde FV_F = F$) to $V_F(s)$. Since $V_F(s)$ is invertible,
we see that $F(s) = \widetilde FV_F(s)$ is invertible in $\mathcal B$.
Since $F$ satisfies the
condition (I), there is $s':Y \to X$ such that $F(s') = (F(s))^{-1}$. 

Now $F$ is objective, thus the functor $\widetilde F$ is faithful
according to Theorem \ref{1.1}. Since $\widetilde F$
is faithful, it follows from  
$\widetilde FV_F(s') = (\widetilde FV_F(s))^{-1}$ that also 
$V_F(s') = (V_F(s)) ^{-1}.$ 

\subsection{\bf The case when $\mathcal A$ is a Fitting category.}
Let us assume now that $\mathcal A$ is a Fitting category. If $a$ is
an endomorphism of $X = X'\oplus X''$ with $a(X') \subseteq X', a(X'') \subseteq X''$, 
then we write $a = a'\oplus a''$. 

Consider a triangle functor $F$ satisfying the condition (I) and let us write $V = V_F$.
Let $s:X \to Y$
be a morphism in $\mathcal A$ and assume that $V(s)$ is invertible. Thus
also $F(s) = \widetilde FV(s)$ is invertible. Since the condition (I) is satisfied,
there is a morphism $t:Y \to X$ such that $F(s)^{-1} = F(t).$

Let $a = ts$. This is a morphism $X \to X$ and since $\mathcal A$ is
a Fitting category, there is a direct decomposition $X = X'\oplus X''$ with
$a(X') \subseteq X', a(X'') \subseteq X''$ (thus $a = a'\oplus a''$)
such that the restriction $a'$ of $a$ to $X'$
is an automorphism, whereas the restriction $a''$ of $a$ to $X''$ is nilpotent. 
Applying $F$, we see that $1_{F(X)} = F(a) = F(a')\oplus F(a'').$ Since $F(a'')$
is nilpotent, it follows that $F(X'') = 0$ and therefore $V(X'') = 0.$  
We denote by $u':X' \to X$ the
canonical inclusion, by $p':X \to X'$ the canonical projection, thus
$p'u' = 1_{X'}$ and $p'au' = a'$. Now 
$(X',X,X'',u',p'',0)$ and $(X'',X,X',u'',p',0)$ are distinguished triangles in $\mathcal A$.
It follows that $(V(X'),V(X),0,V(u'),0,0)$ and $(0,V(X),V(X'),0,V(p'),0)$ are distinguished
triangles in $\mathcal A/\Ker F$. It follows that $V(p')$ and $V(u')$ are
isomorphisms, thus $V(p')V(u') = 1_{V(X')}$ implies that $V(u')V(p') = 1_{V(X)}.$

Let $b' = (a')^{-1}:X' \to X'.$  We consider the map 
$$
 u'b'p'tsu'p' = u'b'p'au'p' = u'b'a'p' = u'p':X \to X.
$$
If we apply $V$, we get 
$$
 V(u'b'p't)V(s) = 
 V(u'b'p'ts) = V(u'b'p'ts)V(u'p') = V(u'b'p'tsu'p') = V(u'p') = 1_{V(X)}.
$$
This shows that $V(s)^{-1} = V(u'b'p't).$ 

\subsection{Examples of Fitting categories.}\label{Fitting}
Let $k$ be a field. A
{\it $\Hom$-finite $k$-category}
$\mathcal A$ is an additive category such that $\Hom_{\mathcal A}(X,
Y)$ is a finite-dimensional $k$-space for arbitrary objects $X$ and
$Y$ of $\mathcal A$, such that the composition of morphisms is
$k$-bilinear. A $\Hom$-finite $k$-category $\mathcal A$ is called {\it a
Krull-Remak-Schmidt category} provided all idempotents of $\mathcal A$
split. It is well-known that {\it any $\Hom$-finite Krull-Remak-Schmidt
$k$-category is a Fitting category.}

Let us outline the proof. If $\Lambda$ is a finite-dimensional $k$-algebra, 
the classical
Fitting lemma asserts that the category $\mo \Lambda$
of all finite-dimensional $\Lambda$-modules is a Fitting category. Of course,
also the full subcategory $\proj \Lambda$ of all projective $\Lambda$-modules
is a Fitting category. 
Now assume that $\mathcal A$ is a  $\Hom$-finite Krull-Remak-Schmidt
$k$-category. Let $X$ be an object in $\mathcal A$ and $\add(X)$ the
full subcategory of all direct summands of finite direct sums of copies of $X$.
Let $\Gamma(X) = \End(X)^\op.$ Then the category $\add(X)$ is equivalent
to the category $\proj \Gamma(X)$ of all projective $\Gamma(X)$-modules, thus
it is a Fitting category.

\section{\bf Proof of Theorem \ref{1.3}}

(i) $\implies$ (ii). According to Proposition \ref{WSMISM}, the condition
(SM) is equivalent to the conditions (WSM) and (I). According to Theorem
\ref{1.1}, the condition (WSM) implies that $F$ is objective. Thus $F$ is
an objective functor which satisfies the condition (I). According 
to Theorem \ref{1.2}(2), we see that $V_F$ is full.
   \medskip

(ii) $\implies$ (iii).  We assume that $F$ is objective and $V_F$ is full. 
We always have $\Ker(F) = \Ker(V_F)$. Lemma \ref{objective} asserts that 
$\ker(F) = \langle \Ker(F)\rangle,$ since by assumption $F$ is objective. 
Since $V_F$ is always objective, we similarly have $\ker(V_F) = \langle \Ker(V_F)\rangle.$
Thus, we see that $\ker(F) = \ker(V_F)$. It follows that there exists a faithful
functor $\Phi:\mathcal A/\ker(F) \to \mathcal A/\Ker(F)$ such that
$V_F = \Phi\pi_F$ (namely $\Phi(\overline f) = V_F(f)$, where $\overline f$
is the residue class of $f$ modulo $\ker(F)$).

Since $V_F$ is full, the factorization 
$V_F = \Phi\pi_F$ shows that also $\Phi$ is full. Altogether we see that 
$\Phi$ is full and faithful. Of course, $\Phi$ is dense,
since the objects of both $\mathcal A/\ker(F)$ and $\mathcal A/\Ker(F)$ are those
of $\mathcal A$ and are not permuted under $\Phi$. 

   \medskip 
(iii) $\implies$ (ii). Let $\Phi:\mathcal A/\ker(F) \to \mathcal A/\Ker(F)$
be an equivalence with $V_F = \Phi\pi_F.$ 
Since $\pi_F$ is always full, and $\Phi$ is an equivalence of functors, also
$V_F$ is full.

Since $V_F$ is always objective, and $\Phi^{-1}$ is an equivalence,
also $\pi_F = \Phi^{-1}V_F$ is objective. But this implies that $F$ is objective.
(Namely, assume that $F(f) = 0$, then $f \in \ker(F) = \ker(\pi_F)$. Since
$\pi_F$ is objective, $f$ factors through an object $K$ with $\pi_F(K) = 0.$
Thus $1_K$ belongs to $\ker(F)$, but this means that $F(K) = 0.$) 
     \medskip  

(ii) $\implies$ (iv). We assume that $F$ is objective and that $V_F$ is full.
There is the factorization $F = \widetilde FV_F$. Let $F_1 = V_F$ and
$F_2 = \widetilde F$. By assumption, $F_1 = V_F$ is a full and dense
triangle functor. Since $F$ is objective,
Theorem \ref{1.1} asserts that $F_2 = \widetilde F$ is faithful.
	\medskip 

(v) $\implies$ (i).
We assume that $F = F_2F_1:\mathcal A \to \mathcal B$ 
where $F_1$ is a full triangle functor
whereas $F_2$ is a faithful triangle functor.
In order to show (SM), we start with a morphism $u:X \to Y$ in $\mathcal A$ such that
$F(u)$ is a splitting monomorphism. Let $(X,Y,Z,u,v,w)$ be a distinguished 
triangle in $\mathcal A,$ thus \newline 
$(F(X),F(Y),F(Z),F(u),F(v),\xi_XF(w))$ is a 
distinguished triangle in $\mathcal B$. Since $F(u)$ is a splitting monomorphism,
we know that $\xi_XF(w) = 0$, thus also $F(w) = 0.$  Since $F_2$ is
faithful, it follows from $F_2F_1(w) = 0$ that $F_1(w) = 0$ and therefore
$F_1(u)$ is a splitting monomorphism. In this way, we see that there is 
a morphism $c:F_1(Y) \to F_1(X)$ such that $cF_1(u) = 1_{F_1(X)}$. Since
$F_1$ is full, there is $u':Y \to X$ such that $F_1(u') = c,$ thus
$F_1(u'u) = 1_{F_1(X)}.$ We apply $F_2$ to this equality in order to see that
$F(u'u) = 1_{F(X)}.$  $\s$
	
%%%%%%%%%%%%%%%%%%%%%%%%%%%%%%%%%%%%%%%%%%%%%%%%%%%%%%%%%%%%%%%%%%%%%%%%%%%%%%%%
\section{\bf Dual conditions}

\subsection{}
We also may consider dual conditions, 
in particular the following ones:

(WSE) For each morphism $v: Y \longrightarrow Z$ in
$\mathcal A$ such that $F(v)$ is a splitting epimorphism in
$\mathcal B$, there exists a morphism $v': Z' \longrightarrow Y$
such that $F(vv')$ is an isomorphism in $\mathcal B$.

(SE) For each morphism $v: Y \longrightarrow Z$ in
$\mathcal A$ such that $F(v)$ is a splitting epimorphism in
$\mathcal B$, there exists a morphism $v': Z \longrightarrow Y$
such that $F(vv') = 1_{F(Z)}$.

(RSE) For each morphism $v: Y \longrightarrow Z$ in
$\mathcal A$ such that $F(v)$ is a splitting epimorphism in
$\mathcal B$, the morphism $v$ is a splitting  epimorphism in $\mathcal A$.
	  \medskip

Of course, there are the following trivial implications:
$$
 \text{(RSE)} \implies \text{(SE)} \implies \text{(WSE)}
$$
Note that most of the conditions considered in the paper are self-dual conditions:
that a functor $F$ is objective or that $V_F$ is full, or that
$F$ is faithful, are self-dual conditions. 
Here, ''duality'' (or better: left-right symmetry) refers to 
the procedure of looking at the opposite $\mathcal A^\op$ of a given category
$\mathcal A$, and to consider a functor $F:\mathcal A \to \mathcal B$ as 
a functor $F^\op:\mathcal A^\op \to \mathcal B^\op$, with $F^\op(X) = F(X), F^\op(F) = F(f)$
for any object $X$ and any morphism $f$ in $\mathcal A^\op$. Since we assume
that $F$ is covariant, also $F^\op$ is covariant. 
For example, we see that $F$ satisfies (SE) 
if and only if it satisfies both (WSE) and (I), this is the dual assertion of Proposition
\ref{WSMISM}.
	\medskip

By duality, the theorems \ref{1.1}, \ref{1.3} and the proposition \ref{faithful-RSM}
yield:
	\medskip

%%%%%%%%%%%%%%%%%%%%%%Theorem 1.1'%%%%%%%%%%%%%%%%%%%%%%%%%%%%%
{\bf Theorem 1.1$'$.}  Let $F: \mathcal
A \longrightarrow \mathcal B$ be a triangle functor between
triangulated categories. Then the following are equivalent:

${\rm (i)}^\op$ \ \ $F$ satisfies the condition ${\rm (WSE)}$.

$\rm (ii)$ \ \ $F$ is objective;
     \medskip

     %%%%%%%%%%%%%%%%%%%%%%Theorem 1.3'%%%%%%%%%%%%%%%%%%%%%%%%%%%%%
     {\bf Theorem 1.3$'$.}   Let $F: \mathcal A \longrightarrow \mathcal B$ be a
     triangle functor between triangulated categories. Then the following
     are equivalent:

     ${\rm (i)}^\op$ \ \  $F$ satisfies the condition ${\rm (SE)}$;

     $\rm (ii)$ \ \ $F$ is objective and $V_F$ is full.
     \medskip

{\bf Proposition 3.2$'$.}  Let $F: \mathcal
A \longrightarrow \mathcal B$ be a triangle functor between
triangulated categories. Then the following are equivalent:

${\rm (i)}^\op$ \ \ $F$ satisfies the condition ${\rm (RSE)}$.

$\rm (ii)$ \ \ $F$ is faithful.

     \medskip

In particular, we see that $F$ satisfies the condition (WSM) if and only if
it satisfies the condition (WSE), that $F$ satisfies the condition (SM) if and only if
it satisfies the condition (SE) and that $F$ satisfies the condition (RSM) if and only if
it satisfies the condition (RSE).
   \medskip

\subsection{}
As a bonus for the reader, let us insert a direct proof that the condition (WSM) for a
triangle functor $F$ 
implies the condition (WSE):
	\medskip 

\noindent {\bf Proof.} \ Assume that $F$ is a triangle functor which satisfies the condition
(WSM). Given a morphism $v: Y\longrightarrow Z$ such that
$F(v)$ is a splitting epimorphism, consider a distinguished
$(X,Y,Z,u,v,w).$ Applying $F$ we know that $F(v)$ is a splitting epimorphism, thus
$F(u)$ is a splitting
monomorphism. By (WSM) there exists a morphism $u': Y
\longrightarrow X'$ such that $F(u')F(u)$ is an isomorphism in
$\mathcal B$. We embed $u'$ into a distinguished triangle 
$(Y,X',Z'[1],u',w',v'[1])$. By the octahedral axiom we
get the following commutative diagram
$$
\xymatrix{X \ar[r]^-u \ar @{=}[d] &Y
\ar[d]^-{u'}\ar[r]^-v &Z\ar@{-->}[d] \ar[r]^-w &X[1] \ar @{=}[d]
\\X \ar[r]^-{u'u} &X' \ar[d]^-{w'} \ar[r]&K\ar@{-->}[d] \ar[r] &X[1]\ar
[d]^{u[1]}\\  & Z'[1] \ar@{=}[r]\ar[d]^-{v'[1]} &Z'[1] \ar[r]^-{v'[1]}\ar[d]^{\beta[1]}  &Y[1]\\
& Y[1] \ar[r]^-{v[1]} &Z[1] & & .}$$ Since $F(u'u)$ is an isomorphism, 
it follows that $F(K) = 0,$ and hence
$F(\beta[1])$ is an isomorphism. Thus $F(vv') = F(v)F(v') =
F(\beta)$ is an isomorphism. This shows that $F$ satisfies the condition (WSE).
 $\s$

 %%%%%%%%%%%%%%%%%%%%%%%%%%%%%%%%%%%%%%%%%%%%%%%%%%%%%%%%%%%%%%
 \section{\bf Examples of triangle functors which
 are not objective}\label{sincere-notobjective}

 The aim of this section is to present examples of 
triangle functors which are not objective (but sincere).
	 \medskip 

If one compares triangulated categories with abelian categories,
then one relates the triangle functors between triangulated
categories to the exact functors between abelian categories,
these are the functors which preserve the given structure. 
Whereas there do exist triangle functors which are not objective,
all exact functors are objective, as the following lemma shows.

\begin{lem}\label{exactobj} Let $F: \mathcal A \longrightarrow \mathcal B$ be an exact
functor between abelian categories. Then $F$ is objective. Thus, an
exact sincere functor between abelian categories is faithful.
\end{lem}

\noindent{\bf Proof.} Let $f: X \longrightarrow Y$ be a morphism in
$\mathcal A$ such that $F(f) = 0.$ Let $I$ be the image of $f$, say
$f = hg$ with $g: X \longrightarrow I$ an epimorphism and $h: I
\longrightarrow Y$ a monomorphism. Then $F(f) = F(hg) = F(h)F(g)$.
Since $F$ is exact, $F(g)$ is epic and $F(h)$ is monic. Thus $F(I)$
is the image of $F(f)$. Since $F(f) = 0$, it follows that $F(I) =
0.$ By definition $F$ is objective. $\s$

\begin{prop}\label{nonobjective} Let $F_0: \mathcal A
\longrightarrow \mathcal B$ be an exact sincere functor between
abelian categories.  Then $F_0$ induces a sincere triangle functor
$F: D^b(\mathcal A) \longrightarrow  D^b(\mathcal B),$
where $D^b(\mathcal A)$ is the bounded derived category of $\mathcal
A$. Moreover, if $\mathcal A$ is not semi-simple whereas $\mathcal B$ is
semi-simple, then  $F$ is not objective.
\end{prop}

Let us add that such a functor $F$ always satisfies the condition (I).
Namely, since $F$ is sincere, the Verdier quotient functor $V_F$ 
is the identity functor,
in particular $V_F$ is full. Thus, according to Theorem \ref{1.2}, $F$ satisfies
the condition (I).
    \medskip

\noindent{\bf Proof.}  Since $F_0: \mathcal A \longrightarrow
\mathcal B$ is an exact functor between abelian categories, it
induces a triangle functor $F = (F, {\rm Id}): D^b(\mathcal A)
\longrightarrow  D^b(\mathcal B)$, which maps a complex $C$ with
cohomology ${\rm H}^n(C)$ to the complex $F(C)$ with cohomology
$F_0({\rm H}^n(C)) = {\rm H}^n(F(C)).$

Assume that $F(C)  = 0.$ Then $F_0({\rm H}^n(C)) = 0,$ for all
$n\in\Bbb Z$. Since $F_0$ is sincere, it follows that ${\rm H}^n(C) =
0,$ for all $n\in\Bbb Z$, thus $C$ is acyclic and therefore $C =
0$ in $D^b(\mathcal A).$ This shows that $F$ is sincere.

Since $\mathcal A$ is not semi-simple, there exist object $X$ and
$Y$ in $\mathcal A$ with $\Ext^1_\mathcal A(X, Y) \neq 0$. Since
$\mathcal B$ is semi-simple, $\Ext^1_\mathcal B(F(X), F(Y)) = 0.$
Thus $\Hom_{D^b(\mathcal A)}(X, Y[1]) \neq 0$, but
$\Hom_{D^b(\mathcal B)}(F(X), F(Y)[1]) = 0.$ That is, $F$ is not
faithful. It follows that $F$ cannot be objective, since sincere
objective functors are faithful. $\s$
	  \medskip

{\bf Example.} 
Let us consider an example in detail. Let $A$ be the path
algebra of the quiver $b \longrightarrow a$ over the field $k$ and let
$B$ be the semisimple algebra given by the quiver with the two
vertices $a, b$ and no arrow. Note that $B$ is a subalgebra of $A$ and we
consider the forgetful functor $F_0:
A\mbox{-}{\rm mod} \longrightarrow B\mbox{-}{\rm mod}$, given by the
inclusion map $B \to A$. 

Given a vertex $x$, we denote by $S_A(x)$ or $S_B(x)$ the simple $A$-module
of $B$-module, respectively, corresponding
to the vertex $x$, and we denote by $P_A(x)$ the indecomposable projective $A$-module 
corresponding to the vertex $x$. The functor $F_0$ sends
$S_A(x)$ to $S_B(x)$ for $x = a, b$, and it sends $P_A(b)$ to $S_B(a)\oplus
S_B(b)$. Clearly, $F_0$ is an exact and faithful functor.

The upper part of the following picture shows the Auslander-Reiten
quiver of $D^b(A)$, the dashed lines indicate the mesh relations.
The lower part is the Auslander-Reiten quiver of $D^b(B)$ (it just consists
of isolated vertices) and here we use
dotted lines to indicate the two shift orbits in $D^b(B)$ (in the upper part, the
shift orbits are not marked in this way).

\[\xymatrix@C=0.2pc @R=2pc {\quad\quad\ar[dr] \ar@{--}[rr]&  & \;\;P_A(b)\;\; \ar[dr]^-{v}\ar@{--}[rr] & &
S_A(a)[1]\ar[dr]\ar@{--}[rr] & & S_A(b)[1]\ar[dr]\ar@{--}[rr]& & \\
\ar@{--}[r] & S_A(a)\ar[ur]^-{u}\ar@{--}[rr]&&
S_A(b)\ar[ur]^-{w}\ar@{--}[rr]& & P_A(b)[1]\ar[ur]\ar@{--}[rr]& &
S_A(a)[2]\ar[ur]\quad\quad \\ \\
\ar@{..}[r] & \;S_B(a)\;\ar@{..}[rrr]&& & S_B(a)[1] \ar@{..}[rrr]& & &
S(a)[2]\\
\ar@{..}[rrr] & & & \;S_B(b)\;\ar@{..}[rrr]& & & S_B(b)[1] \ar@{..}[rr]& &
}\]

\noindent
The induced functor $F: D^b(A) \longrightarrow D^b(B)$ sends
$S_A(x)[i]$ to $S_B(x)[i]$ for $x = a,b$ and all $i\in \Bbb Z$, and it
sends $P_A(b)[i]$ to $S_A(a)[i]\oplus S_A(b)[i].$ In $D^b(A)$ we have
labeled three arrows $u, v, w$, they form a distinguished triangle
$(S(a), P(b), S(b), u, v, w)$.
	\medskip

	Consider the map $w: S_A(b) \longrightarrow S_A(a)[1].$ Since
	$\Hom_{D^b(B)}(S_B(b), S_B(a)[1]) = 0$, we have $F(w) = 0$. Thus, we see
	that $F$ is not faithful.

	On the other hand, consider the map $u: S_A(a) \longrightarrow P_A(b)$.
	Applying the functor $F$, we obtain the inclusion map $S_B(a)
	\longrightarrow S_B(a)\oplus S_B(b)$ which is splitting mono: there is a
	projection map $u': S_B(a)\oplus S_B(b) \longrightarrow S_B(a)$ with $u'
	F(u) = 1_{S_B(a)}$. Since there is no non-zero map $P_A(b)
	\longrightarrow S_A(a)$, such a map $u'$ is not in the image of $F$.
	This shows that the condition (SM) is not satisfied. Thus, Theorem \ref{1.1}
asserts that $F$ is not objective.
 \newpage 

\section{\bf Overlook}

\subsection{The main conditions.}

We consider any triangle functor $F:\mathcal A \to \mathcal B$. 
	\medskip

\includegraphics{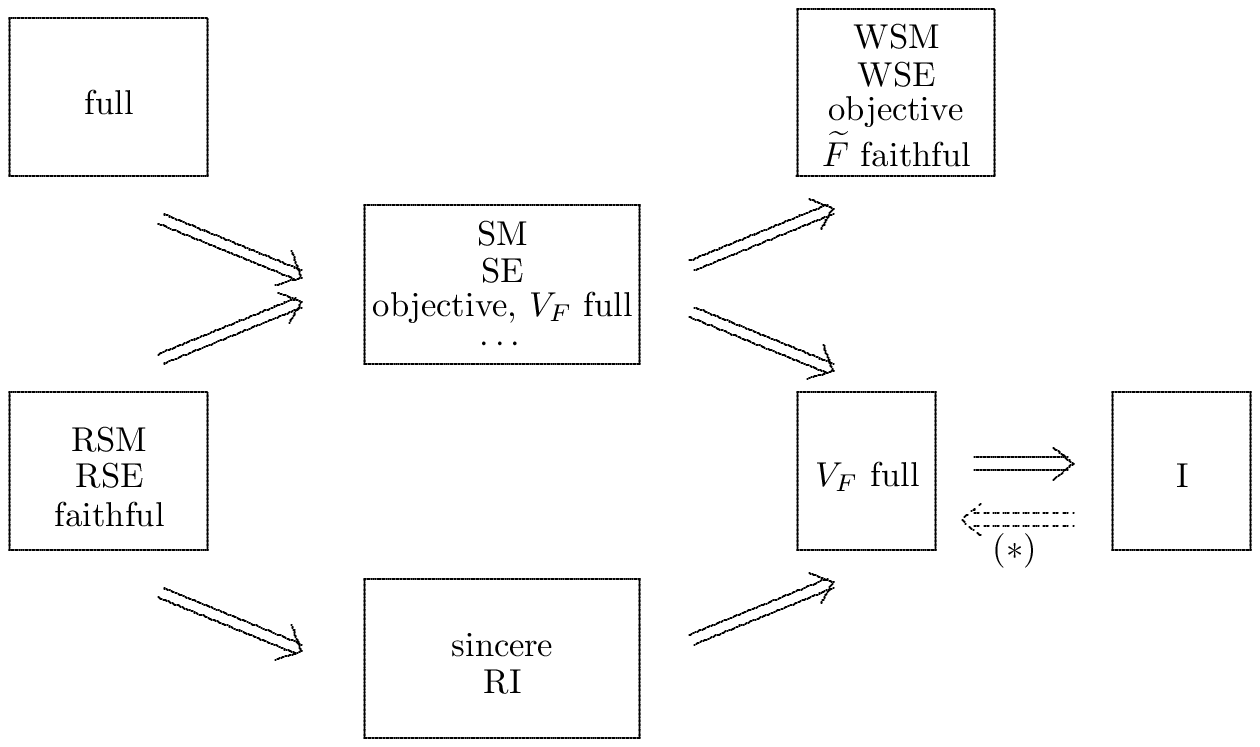}
	\medskip 

\noindent 
The conditions in any box, line by line, are equivalent; note that in 
the central box which mentions the conditions (SM)
and (SE), the dots indicate that there are several further equivalent
conditions, namely the conditions (iii), (iv) and (v)
mentioned in Theorem  \ref{1.2} as well as the conjunction of the conditions
(WSM) and (I), see Proposition \ref{WSMISM}, and dually also the conjunction of
(WEM) and (I). 

The arrows show the relevant implications between the boxes.
The dashed implication with the label $(*)$ 
is valid under the assumption that $F$ is objective or that 
$\mathcal A$ is a Fitting category.

\newpage

\subsection{References.} Two implications are trivial: a faithful functor is of course sincere. And if $F$
is a sincere triangle functor, then $V_F$ is the identity functor, thus full.
Here are the references for the remaining implications mentioned above:
	\medskip

\includegraphics{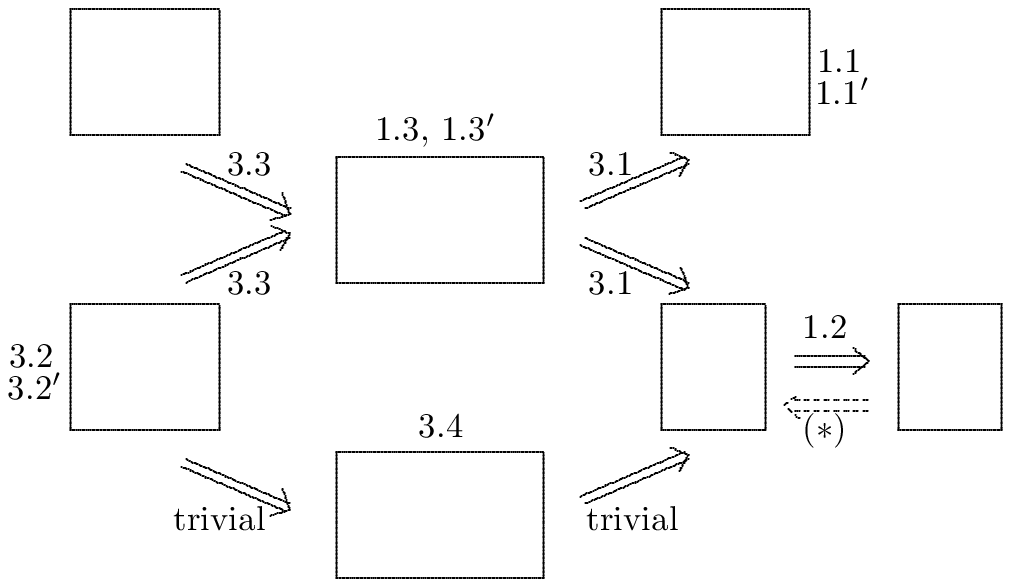}

\subsection{Examples.} Finally, let us outline typical examples in order to see
that the implications (A), (B), (C), (D), (E) and (F) cannot be reversed:
	\medskip
     
\includegraphics{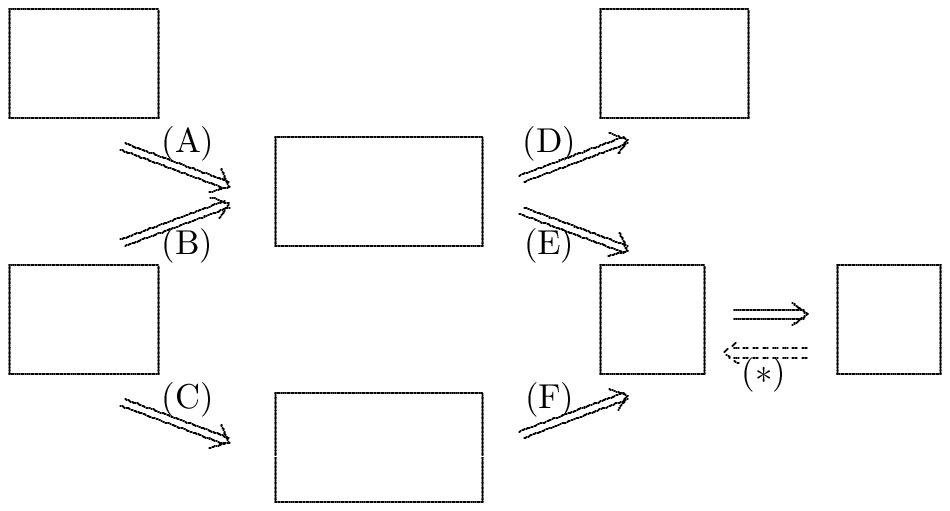}
	\medskip

(A) Take any faithful functor which is not full, for example the zero functor
$0 \to \mathcal A$, where $\mathcal A$ is a non-zero triangulated category. 

(B) Take any full functor which is not faithful, for example the zero functor
$\mathcal A \to 0$, where $\mathcal A$ is a non-zero triangulated category. 

(C) Take any sincere functor $F$ which is not objective as presented in section 
\ref{sincere-notobjective}. Such a functor is of course not faithful.

(D) In order to find an objective functor $F$ such that $V_F$ is not full, 
consider a Verdier quotient functor $V_{\mathcal K}$, these functors are
very seldom full! 

For example, let $A$ be an Artin algebra, $\mo A$ the
category of finitely generated $A$-modules, $K^-(A)$ the homotopy
category of the upper bounded complexes over $\mo A$, $\mathcal E$
the full subcategory of $K^-(A)$ consisting of the upper bounded
acyclic complexes, and $D^-(A)$ the derived category of the upper
bounded complexes over $\mo A$. Then we have the Verdier quotient functor
$V_\mathcal E: K^-(A) \longrightarrow K^-(A)/\mathcal E
= D^-(A)$. It is well-known that $V_\mathcal E$ is full if and only
if $A$ is semi-simple.

(E) Take any sincere functor $F$ which is not objective as presented in section 
\ref{sincere-notobjective}. Since $F$ is sincere, it satisfies the condition (I)
but it cannot satisfy the condition (SM), since otherwise it would be objective.

(F) Take a functor $F$ which is not sincere, such that $V_F$ is full, for
example  the zero functor
$\mathcal A \to 0$, where $\mathcal A$ is a non-zero triangulated category.

   \vskip20pt

    {C. M. Ringel \par
                Department of Mathematics \par
		            Shanghai Jiao Tong University \par
			                Shanghai 200240, P. R. China, \par
					            and \par
						                King Abdulaziz University, P O Box 80200\par
								            Jeddah, Saudi Arabia\par
									                e-mail: {ringel$\symbol{64}$math.uni-bielefeld.de} \par
											            \medskip
												                P. Zhang\par
														            Department of Mathematics \par
															                Shanghai Jiao Tong University \par
																	            Shanghai 200240, P. R. China, \par
																		                e-mail: {pzhang$\symbol{64}$sjtu.edu.cn} }
																				\end{document}